\documentclass[11pt]{article}

\usepackage{latexsym,color,graphicx,amsthm,amsmath,amssymb}
\usepackage{lineno}

\def\reals{{\mathbb R}}

\newtheorem{theorem}{Theorem}[section]
\newtheorem{corollary}[theorem]{Corollary}

\newtheorem{lemma}[theorem]{Lemma}

\newtheorem{question}[theorem]{Question}

\makeatletter
\newcommand{\ProofEndBox}{{\ifhmode\unskip\nobreak\hfil\penalty50 \else
          \leavevmode\fi\quad\vadjust{}\nobreak\hfill$\Box$
            \finalhyphendemerits=0 \par}}
\makeatother

\newcommand{\proofend}{\ProofEndBox\smallskip}
\newcommand{\ignore}[1]{}
\begin{document}


\title{Highly incidental patterns on a quadratic hypersurface in $\mathbb{R}^4$\thanks{%
Work on this paper by Noam Solomon was supported by Grant 892/13
from the Israel Science Foundation. Ruixiang Zhang was supported by Princeton University. Part of this research was
performed while the authors were visiting the Institute for Pure and
Applied Mathematics (IPAM), which is supported by the National
Science Foundation. }}

\author{
Noam Solomon\thanks{%
School of Computer Science, Tel Aviv University, Tel Aviv 69978,
Israel. {\sl noam.solom@gmail.com} } \and
Ruixiang Zhang\thanks{%
Department of Mathematics, Princeton University, Princeton, NJ 08540
{\sl ruixiang@math.princeton.edu} } }

\maketitle

\begin{abstract}
In~\cite{SS4d}, Sharir and Solomon showed that the number of
incidences between $m$ distinct points and $n$ distinct lines in
$\reals^4$ is \begin {equation}\label{fourdimbound}
O^*\left(m^{2/5}n^{4/5}+ m^{1/2}n^{1/2}q^{1/4} +
m^{2/3}n^{1/3}s^{1/3} + m + n\right), \end {equation} provided that
no 2-flat contains more than $s$ lines, and no hyperplane or
quadric contains more than $q$ lines, where the $O^*$ hides a
multiplicative factor of $2^{c\sqrt {\log m}}$ for some absolute
constant $c$.

In this paper we prove that, for integers $m, n$ satisfying $n^{9/8}
< m < n^{3/2}$, there exist $m$ points and $n$ lines on the
quadratic hypersurface in $\mathbb{R}^4$
\begin{equation*}
\{(x_1,x_2,x_3,x_4)\in \reals^4 \mid x_1 = x_2^2 + x_3^2 - x_4^2\},
\end{equation*} such that (i) at most $s=O(1)$ lines lie on any
2-flat, (ii) at most $q=O(n/m^{1/3})$ lines lie on any hyperplane,
and (iii) the number of incidences between the points and the lines
is $\Theta(m^{2/3}n^{1/2})$, which is asymptotically larger than the
upper bound in~(\ref{fourdimbound}), when $n^{9/8}<m<n^{3/2}$. This
shows that the assumption that no quadric contains more than $q$
lines (in the above mentioned theorem of~\cite{SS4d}) is necessary
in this regime of $m$ and $n$.

By a suitable projection from this quadratic hypersurface onto
$\reals^3$, we obtain $m$ points and $n$ lines in $\reals^3$, with
at most $s=O(1)$ lines on a common plane, such that the number of
incidences between the $m$ points and the $n$ lines is
$\Theta(m^{2/3}n^{1/2})$. It remains an interesting question to
determine if this bound is also tight in general.
\end{abstract}


\section{Introduction}

Let $P$ be a set of $m$ distinct points in $\reals^2$ and let $L$ be
a set of $n$ distinct lines in $\reals^2$. Let $I(P,L)$ denote the
number of incidences between the points of $P$ and the lines of $L$;
that is, the number of pairs $(p,\ell)$, such that $p\in P$,
$\ell\in L$ and $p\in \ell$. The classical
Szemer\'edi--Trotter theorem~\cite{SzT} yields the worst-case tight
bound

\begin{equation} \label{inc2}
I(P,L) = O\left(m^{2/3}n^{2/3} + m + n \right) .
\end{equation}
This bound clearly also holds in three, four, or any higher
dimensions which can be easily proved by projecting the given lines and points onto some
generic plane. Moreover, the bound will continue to be worst-case
tight by placing all the points and lines in a common plane, in a
configuration that yields the planar lower bound.

In the groundbreaking paper of Guth and Katz~\cite{GK2}, an improved
bound has been derived for $I(P,L)$, for a set $P$ of $m$ points and
a set $L$ of $n$ lines in $\reals^3$, provided that not too many
lines of $L$ lie in a common plane \footnote{The additional
requirement in \cite{GK2}, that no regulus contains too many lines,
is not needed for the bound given below.}. Specifically, they
showed:
\begin{theorem}  [Guth and Katz~\protect{\cite{GK2}}]
\label {ttt} Let $P$ be a set of $m$ distinct points and $L$ a set
of $n$ distinct lines in $\reals^3$, and let $s\le n$ be a
parameter, such that
no plane contains more than $s$ lines of $L$. Then
$$
I(P,L) = O\left(m^{1/2}n^{3/4} + m^{2/3}n^{1/3}s^{1/3} + m +
n\right).
$$
\end{theorem}
\noindent{\bf Remark.} When $s=\Theta(\sqrt n)$, this bound is known
to be tight, by a generalization to three dimensions of Elekes'
planar construction of points and lines on an integer grid (see Guth
and Katz~\cite{GK2} for the details). For smaller values of $s$, it
is an open problem to give lower bounds or improve the upper bound,
and the case $s=O(1)$ is of particular interest. In
Theorem~\ref{th:mainthree} we give an improved upper bound, and it
remains a question (see Question~\ref{conj}) whether it is tight.

In a recent paper of Sharir and Solomon~\cite{SS4d}, the following
analogous and sharper result in four dimensions was established.

\begin{theorem} \label{th:main0}
Let $P$ be a set of $m$ distinct points and $L$ a set of $n$
distinct lines in $\reals^4$, and let $q,s\le n$ be parameters, such
that (i) each hyperplane or quadric contains at most $q$ lines of
$L$, and (ii) each 2-flat contains at most $s$ lines of $L$. Then
\begin {equation}
\label {th:shso} I(P,L) \le 2^{c\sqrt{\log m}} \left( m^{2/5}n^{4/5}
+ m \right) + A\left( m^{1/2}n^{1/2}q^{1/4} + m^{2/3}n^{1/3}s^{1/3}
+ n\right) ,
\end {equation}
where $A$ and $c$ are suitable absolute constants. When $m\le
n^{6/7}$ or $m\ge n^{5/3}$, there is the sharper bound
\begin {equation}
\label {ma:in0x} I(P,L) \le A \left( m^{2/5}n^{4/5} + m +
m^{1/2}n^{1/2}q^{1/4} + m^{2/3}n^{1/3}s^{1/3} + n\right) .
\end {equation}
In general, except for the factor $2^{c\sqrt{\log m}}$, the bound is
tight in the worst case, for any values of $m,n$, and for
corresponding suitable ranges of $q$ and $s$.
\end{theorem}
The term $m^{2/3}n^{1/3}s^{1/3}$ comes from the planar
Szemer\'edi--Trotter bound (\ref{inc2}), and is unavoidable, as it
can be attained if we densely pack points and lines into 2-flats, in
patterns that realize the bound in (\ref{inc2}).

Likewise, the term $m^{1/2}n^{1/2}q^{1/4}$ comes from the bound of
Guth and Katz~\cite{GK2} in three dimensions (as in
Theorem~\ref{ttt}), and is again unavoidable, as it can be attained
if we densely ``pack'' points and lines into hyperplanes, in
patterns that realize the bound in three dimensions.

In this paper we show that the condition in assumption (i) of
Theorem~\ref{th:main0} that quadrics also do not contain too many
lines, cannot be dropped, by proving the following theorem.

\begin{theorem}
\label{th:main} For each positive integer $k$ and each $\alpha > 0$,
there exists $m = \Theta(k^{3+3 \alpha})$ points and $n = \Theta(k^{2
+ 4 \alpha})$ lines on the quadratic hypersurface
\begin{equation*}
S:=\{(x_1,x_2,x_3,x_4)\in \reals^4 \mid x_1 = x_2^2 + x_3^2 -
x_4^2\}
\end{equation*}
in $\mathbb{R}^4$, such that there are at most $O(1)$ lines lying on
any 2-flat and $O(k^{1 + 3 \alpha})$ lines lying on any hyperplane,
and $I(P,L)=\Theta(k^{3 + 4 \alpha})$. \\
\end{theorem}

Given integers $m$ and $n$, there are $k,\alpha$ such that $m =
\Theta(k^{3+3 \alpha})$ points and $n = \Theta(k^{2 + 4 \alpha})$.
Substituting these values in Theorem~\ref{th:main}, we obtain the
following corollary.

\begin {corollary}
\label{co:main} For integers $m,n$, there is a configuration of $m$
points and $n$ lines in $\reals^4$, such that all the points (resp.,
lines) are contained (resp., fully contained) in $S$, and (i) the
number of lines in any common 2-flat is $O(1)$, (ii) the number of
lines in a common hyperplane is $O(n/m^{1/3})$, and (iii) the number
of incidences between the points and lines is
$\Omega(m^{2/3}n^{1/2}+m+n)$.
\end {corollary}

\noindent{\bf Remarks. (1)} For integers $m,n,$ satisfying
$n^{9/8}<m<n^{3/2}$, the number incidences $\Omega(m^{2/3}n^{1/2})$
in Corollary~\ref{co:main} is asymptotically larger than the bound
of Theorem~\ref{th:shso} for the number of incidences
$O(m^{2/5}n^{4/5}+m^{1/2}n^{1/2}q^{1/4}+m^{2/3}n^{1/3}+m+n) =
O(m^{2/5}n^{4/5}+m^{5/12}n^{3/4} + m + n)$ (as $q=O(n/m^{1/3})$).
This implies that the condition in assumption (i) of
Theorem~\ref{th:main0} cannot be dropped, in this regime of $m$ and
$n$.

\smallskip

\noindent {\bf (2)} We note that the number of 2-rich points
determined by $n$ lines in $\reals^4$ is $O(n^{3/2})$, provided that
at most $O(\sqrt n)$ of the lines lie on a common plane or
\emph{regulus}, \footnote{A regulus is a quadratic surface that is
doubly ruled by lines. For more details about reguli, see e.g.,
Sharir and Solomon~\cite{SS3dv}.} To see this, project the lines
onto some (generic) hyperplane $H$, such that no two lines are
projected onto the same line, and similarly, no two 2-rich points
are projected onto the same 2-rich point, and such that at most
$O(\sqrt n)$ lines lie on a common plane or regulus. Then, the
number of 2-rich points in the configuration of $n$ lines in
$\reals^4$ is equal to the number of 2-rich points in the
configuration of the projected lines onto $H$. By Guth and
Katz~\cite{GK2}, the number of 2-rich points determined by the
projected lines is $O(n^{3/2})$, and therefore the same holds for
the number of 2-rich points in the original configuration of lines
in $\reals^4$. We also notice that in a configuration of $m$ points
and $n$ lines in $\reals^4$, the 1-rich points (i.e., points that
are incident to exactly one line) contribute at most $m$ incidences.
Therefore, in Corollary~\ref{co:main}, as $s=O(1)$, the assumption
that $m \le n^{3/2}$ causes no loss of generality.

\paragraph{Proof Techniques.}
It is a common practice to take geometric objects to be integer
points on certain hypersurfaces (especially quadratic ones) and
varieties passing through a lot of such points, in order to obtain
lower bounds for their incidences. For some most recent applications
of this method, see \cite{Shef} \cite{T2} \cite{Zhang}. In this
paper we obtain our incidence lower incidence bound by taking
integer points and ``low height'' lines on the above hypersurface
$S$.

\paragraph{Projection to $\reals^3$.}
As remarked above, Guth and Katz~\cite{GK2} proved that the number
of incidences between $m$ points and $n$ lines in $\reals^3$ is
$I(P,L) = O\left(m^{1/2}n^{3/4} + m^{2/3}n^{1/3}s^{1/3} + m +
n\right)$, provided that no plane contains more than $s$ lines of
$L$. When $s=\Theta(\sqrt n)$, this bound is tight, by a
generalization to three dimensions of Elekes' construction of points
and lines on an integer grid in the plane (see Guth and
Katz~\cite{GK2} for the details). For smaller values of $s$, it is
an open problem to give lower bounds or improve the upper bound,
where the case $s=O(1)$, is of particular interest.

By choosing a generic projection from $\reals^4$ to $\reals^3$, we
show that Corollary~\ref{co:main} directly implies the following
Theorem.
\begin {theorem}
\label{th:mainthree} For integers $m,n$, there is a configuration of
$m$ points and $n$ lines in $\reals^3$, such that (i) the number of
lines in any common plane is $s=O(1)$, and (ii) the number of
incidences between the points and lines is
$\Omega(m^{2/3}n^{1/2}+m+n)$.
\end {theorem}

\noindent{\bf Remark.} When $n^{3/4} \ll m \ll n^{3/2}$, the term
$m^{2/3}n^{1/2}$ dominates over $m$ and $n$, showing that in this
regime of $m$ and $n$, the construction in
Theorem~\ref{th:mainthree}, of $m$ points and $n$ lines with $O(1)$
lines in a common plane, yields a super-linear number of incidences.
As observed above, the bound of Guth and Katz~\cite{GK2} implies
that the number of 2-rich points determined by the $n$ lines is
$O(n^{3/2})$, so the assumption that $m \le n^{3/2}$ causes no loss
of generality.

\paragraph{Background.}
Incidence problems have been a major topic in combinatorial and
computational geometry for the past thirty years, starting with the
Szemer\'edi-Trotter bound \cite{SzT} back in 1983. Several
techniques, interesting in their own right, have been developed, or
adapted, for the analysis of incidences, including the
crossing-lemma technique of Sz\'ekely~\cite{Sz}, and the use of
cuttings as a divide-and-conquer mechanism (e.g., see~\cite{CEGSW}).
Connections with range searching and related problems in
computational geometry have also been noted, and studies of the
Kakeya problem (see, e.g., \cite{T}) indicate the connection between
this problem and incidence problems. See Pach and Sharir~\cite{PS}
for a comprehensive survey of the topic.

\ignore{

The landscape of incidence geometry has dramatically changed in the
past six years, due to the infusion, in two groundbreaking papers by
Guth and Katz~\cite{GK,GK2}, of new tools and techniques drawn from
algebraic geometry. Although their two direct goals have been to
obtain a tight upper bound on the number of joints in a set of lines
in three dimensions \cite{GK}, and an almost tight lower bound for
the classical distinct distances problem of Erd{\H o}s \cite{GK2},
the new tools have quickly been recognized as useful for incidence
bounds of various sorts. See \cite{EKS,KMSS,KMS,SSZ,SoTa,Za1,Za2}
for a sample of recent works on incidence problems that use the new
algebraic machinery. }

The simplest instances of incidence problems involve points and
lines. Szemer\'edi and Trotter solved completely this special case
in the plane \cite{SzT}. Guth and Katz's second paper~\cite{GK2}
provides a worst-case tight bound in three dimensions, under the
assumption that no plane contains too many lines; see
Theorem~\ref{ttt}. Under this assumption, the bound in three
dimensions is significantly smaller than the planar bound (unless
one of $m,n$ is significantly smaller than the other), and the
intuition is that this phenomenon should also show up as we move to
higher dimensions. The first attempt in higher dimensions was made
by Sharir and Solomon in ~\cite{surf-socg}. In a recent work, Sharir and Solomon~\cite{SS4d} gave a
tight bound in four-dimensions provided that the number of lines
fully contained in a common hyperplane or quadric is bounded by a
parameter $q$, and the number of lines fully contained in a common
2-flat is bounded by a parameter $s$. Whereas the condition that no
common hyperplane contains more than a bounded number of lines was
known to be necessary, it remained an open question whether the
condition that the number of lines in a common quadric is bounded is
necessary. In this paper, we show that when $n^{9/8}<m<n^{3/2}$,
this condition is indeed necessary, by describing an explicit
quadratic hypersurface in $\reals^4$ containing more incidences than
the bound prescribed by the main theorem of~\cite{SS4d}. This is the
content of Theorem~\ref{th:main}, and Corollary~\ref{co:main}.

We remark that in \cite{T2}, another example of points on a quadratic hypersurface in $\mathbb{F}^4$ with highly incidental pattern was noticed.
There $\mathbb{F}$ is a finite field. Our current quadratic hypersurface and our counting techniques in $\mathbb{R}^4$ are slightly different.
The reader may find it interesting to compare the results here to the results in \cite{T2}.

Another interesting remark is that in three dimensions, there are
certain quadratic surfaces, called reguli, such that if one allows
too many lines to lie on such a regulus, the  number of 2-rich
points determined by them can be larger than the Guth-Katz
bound~\cite{GK2} of $O(n^{3/2})$. The quadratic hypersurface in
$\reals^4$ presented in this paper can be thought of as a higher
degree analogs of regulus. However, If one only cares about
incidences between points and lines (instead of the number of 2-rich
points determined by the lines), the existence of many lines on a
regulus (or any quadratic surface in $\reals^3$) do not yield more
than a linear number of incidences.

\section*{Acknowledgements}
We thank Micha Sharir for his invaluable advice, and the anonymous
referees for their helpful comments.

\section{Proof of Theorem~\ref{th:main}}
\noindent{\bf Proof.} We start by recalling the quadric
\begin{equation}
S=\{(x_1,x_2,x_3,x_4)\in \reals^4 \mid x_1 = x_2^2 + x_3^2 -
x_4^2\},
\end{equation} on which the construction takes place, and define the set of points by
\begin {equation}
P=\{(x_1, x_2, x_3, x_4) \in S \mid x_i \in \mathbb{Z},
i=1,\ldots,4, \quad |x_1| \leq 200 k^{2 + 2 \alpha}, |x_2|, |x_3|,
|x_4| \leq 100 k^{1+ \alpha}\},
\end {equation}
and the set of lines
\begin {equation*}
\begin {array}{ll}
L=\{\{\mathbf{x}+t\mathbf{v}\mid t \in \mathbb{R}\} \subseteq S \mid & \mathbf{x}= (x_1,
\ldots, x_4), \quad \mathbf{v} = (v_1, \ldots , v_4), \quad x_i, v_i \in
\mathbb{Z}, i=1,\ldots, 4, \cr\cr & |x_1| \leq k^{2 + 2 \alpha},
|x_2|, |x_3|, |x_4| \leq k^{1+ \alpha}, \cr\cr & \frac
{k^{1+2\alpha}} 4 \le |v_1| \le 8k^{1 + 2\alpha}, |v_2|, |v_3|\le
k^{\alpha}, \quad v_4^2=v_2^2+v_3^2, \cr\cr &
v_1=2x_2v_2+2x_3v_3-2x_4v_4, \cr\cr & \gcd (v_2, v_3, v_4) = 1,
\text{ and } |v_4| \ge \frac{k^{\alpha}}{2}\}, \end {array}
\end {equation*} for any positive integer $k$ and any $\alpha > 0$.

Since a point on $S$ is uniquely determined by its last three
coordinates, we have $$|P|=|\{(x_2, x_3, x_4) \in \mathbb Z^3 \mid
|x_2|, |x_3|, |x_4| \leq 100 k^{1+ \alpha}\}| = \Theta(k^{3+3\alpha}).$$

The analysis of (an asymptotically tight bound on) the number of
lines of $L$ is a bit more involved. A line
$\{\mathbf{x}+t\mathbf{v}\mid t \in \mathbb{R}\}$ in $L$ (assuming
$\mathbf{x} \in S$, $|x_1| \leq k^{2 + 2 \alpha}, |x_2|, |x_3|,
|x_4| \leq k^{1+ \alpha}$) is fully contained in $S$ if and only if
$$\quad v_1 = 2x_2 v_2 + 2x_3 v_3 - 2x_4 v_4 \text{ and } v_4^2 =
v_2^2 + v_3^2.$$ It follows by Benito and Varona~\cite[Theorem
1]{Be-Va} that the number of primitive integer triples $(v_2, v_3,
v_4)$ (i.e., without a common divisor) satisfying $v_4^2 = v_2^2 +
v_3^2, \quad |v_2|, |v_3| \le k^{\alpha},$ and $|v_4| \ge \frac
{k^{\alpha}} 2$ is $\Theta(k^{\alpha})$. For each such $(v_2, v_3,
v_4)$, we claim that there are $\Omega(k^{3+3\alpha})$ (and
trivially also $O(|P|) = O(k^{3 + 3 \alpha})$) points $\mathbf{x}
\in P$, such that $v_1 = 2x_2 v_2 + 2x_3 v_3 - 2x_4 v_4$ satisfying
$\frac {k^{1+ 2 \alpha}} 4 \le |v_1| \le 8k^{1+ 2 \alpha}$. Indeed,
note that $|v_2|,|v_3| \le |v_4|$. Choosing $|x_2|, |x_3| \le \frac
{|x_4|} 4, \quad \frac {k^{1+\alpha}} 2 \leq |x_4| \leq k^{1+
\alpha}$ (there are at least $\frac {k^{3+3\alpha}} {32}$ choices of
such triples $(x_2,x_3,x_4)$) implies that
$$|2x_2 v_2 + 2x_3 v_3| \le 2|x_2||v_2|+2|x_3||v_3|\le 2\frac {|x_4|} 4(|v_2|+|v_3|)\le |x_4||v_4|,$$ Here
$|v_1|\ge |x_4||v_4|\ge \frac {k^{1+2\alpha}} 4$. The
inequality $|v_1| \le 8k^{1+ 2 \alpha}$ is immediate.

Moreover, each line $\ell$ satisfying the above conditions is
incident to $O(k)$ different points of $P$ (and can thus be
expressed in $O(k)$ different ways as $\{\mathbf{x}+t\mathbf{v}\mid
t \in \mathbb{R}\} \subseteq S,$ for $|x_1| \leq k^{2 + 2 \alpha},
|x_2|, |x_3|, |x_4| \leq k^{1+ \alpha}$). Indeed, parameterize
$\ell$ as $\{\mathbf{x}+t\mathbf{v}\mid t \in \mathbb{R}\}\subset
S$, where $\mathbf{x}, \mathbf{v}$ satisfy
$$|v_1|\ge \frac {k^{1+2\alpha}} 4, \quad |x_1|\le k^{2+2\alpha},$$
and $\mathbf{v}=(v_1,v_2,v_3,v_4)$ is primitive (i.e., its coordinates do not
have a common factor).  Notice that if $|t|>8k$, then the first
coordinate of $\mathbf{x}+t\mathbf{v}$ has absolute value greater than
$k^{2+2\alpha}$, and that if $t\not\in \mathbb Z$, then $\mathbf{x}+t\mathbf{v}\not\in
\mathbb Z^4$ (since $\mathbf{v}$ is primitive and $x \in \mathbb Z^4$). In
either case, $\mathbf{x}+t\mathbf{v} \not\in P$. This implies that
$$\ell\cap P \subseteq \{\mathbf{x}+t\mathbf{v} \mid
 t\in \mathbb Z, \quad |t|\le 8k\},$$ and thus $|\ell\cap P|\le 16k = O(k)$ as claimed.
Therefore, the total number of lines is $ \Omega(\frac{k^{3+ 3
\alpha} k^{\alpha}}{k}) = \Omega(k^{2+ 4 \alpha})$.

It is easy to see that each line in $L$ is incident to $\Omega(k)$ points in $P$. It follows that $|L| = O(k^{2 + 4 \alpha})$. Hence $|L| = \Theta (k^{2 + 4 \alpha})$.

Since each line has $\Theta(k)$ integer points in $P$ on it, we have
$$I(P,L)=\Theta(k^{3+ 4 \alpha}).$$

We now bound the number of lines fully contained in any 2-flat,
 and then bound the number of  lines on any hyperplane. The bounds
will be uniform (i.e., independent of the specific 2-flat or
hyperplane).

Let $\pi$ denote any 2-flat, and we analyze the number of lines that
are fully contained in $\pi \cap S$. We claim that $S$ contains no
planes, so $\pi \not\subset S$. Assume the contrary, then we
parameterize
$$\pi=\{(u_1 s + r_1 t + w_1, u_2 s + r_2 t + w_2, u_3 s + r_3 t +
w_3, u_4 s + r_4 t + w_4) \mid s, t \in \reals\},$$ for constants
$u_i, r_i, w_i \in \reals$, $i = 1, 2, 3, 4$ where $(u_1, u_2, u_3,
u_4)$ and $(r_1, r_2, r_3, r_4)$ are both nonzero and not
proportional to each other. Comparing the coefficients of quadratic
terms in the identity
$$ u_1 s + r_1 t + w_1 \equiv (u_2 s + r_2 t + w_2)^2 + (u_3 s + r_3 t +
w_3)^2 - (u_4 s + r_4 t + w_4)^2, $$ we deduce $(u_2, u_3, u_4)$ and
$(r_2, r_3, r_4)$ are proportional to each other. Hence we may
assume $u_2 = u_3 = u_4 = 0$. But this forces $u_1 = 0$, a
contradiction. Therefore $\pi$ is not contained in $S$. Thus the
intersection $\pi \cap S$ is a curve of degree at most two, so there
are at most two lines fully contained in $\pi \cap S$.

Next, we take any hyperplane $H$, and analyze the number of lines
fully contained in $S\cap H$. The surface $S \cap H$ is a quadratic
2-surface contained in $H$. We will use the classification of (real)
quadratic surfaces in $\reals^3$ (see, e.g., Sylvester's original
paper~\cite{Syl}), and distinguish between two cases.

If the equation of $H$ can be expressed as $x_1 = \varphi(x_2, x_3,
x_4)$, where $\varphi$ is a linear form, then each point $\mathbf{x}
\in H\cap S$ satisfies the equations \ignore{then a line $\ell$ on
is fully contained in $H\cap S$ if and only if every point
$\mathbf{x} \in \ell$ satisfy}
\begin {equation}
\label {eq:quad} \left\{ \begin{array}{l} x_2^2 + x_3^2 - x_4^2 = \varphi(x_2, x_3, x_4),\\
\mathbf{x} \in H. \\
\end{array} \right.
\end {equation}
\ignore{After performing a linear change of coordinates}This is
either a cone, i.e., is linearly equivalent to
$x_2^2+x_3^2-x_4^2=0$, or a hyperboloid of one or two sheets, i.e.,
is linearly equivalent to $x_2^2+x_3^2-x_4^2=1$ or
$x_2^2+x_3^2-x_4^2=-1$, respectively. It is easy to verify (and well
known) that there are no lines on the hyperboloid of two sheets. We
therefore assume that $S \cap H$ is either a cone or a hyperboloid
of one sheet. In these cases, there are at most two lines of $L$
with any given direction that are fully contained in $S \cap H$.
\ignore{Fix any direction $v=(v_2,v_3,v_4)$ (of a line in $S \cap H$
in this parametrization). If the line $\ell=\{\mathbf{x}+t\mathbf{v}
\mid t \in \reals\}$ is fully contained in $S\cap H$, then (under
this linear transform) $(x_2,x_3,x_4)$ satisfy $x_2^2+x_3^2-x_4^2=1$
and $x_2 v_2 + x_3 v_3 - x_4 v_4 = 0$. This is a quadratic curve in
$x_2,x_3,x_4$, and can therefore contain at most two lines (in
$x_2,x_3$ and $x_4$). By applying the inverse of the linear
transform above to return to the original coordinates, it follows
that $S \cap H$ contains at most two lines of the fixed direction.}
Note that if a line $\{\mathbf{x}+t\mathbf{v} \mid t\in\reals \} \in
L$ is fully contained in $S \cap H$, then $v_1 =
\widetilde{\varphi}(v_2, v_3, v_4)$ (where we let
$\widetilde{\varphi}$ denote the linear homogeneous part of
$\varphi$), and $v_4^2 = v_2^2 + v_3^2$ (being the homogeneous part
of degree two in $t$), $|v_2|, |v_3| \le k^{\alpha}$ and $|v_4| \ge
k^{\alpha}$. As observed above, there are $O(k^{\alpha})$ such
triples $(v_2,v_3,v_4)$. Therefore, the number of lines in $L$ that
lie in $S \cap H$ is $O(k^{\alpha})$.

In the remaining case, the equation of $H$ is of the form
$\varphi(x_2, x_3, x_4) = 0$, where $\varphi$ is a linear form. We
can assume, without loss of generality, that the equation of $H$ is
$x_2 = \psi(x_3, x_4)$, where $\psi$ is a linear form (the remaining
case $x_4 = 0$ is simpler to handle). In this case, for every point
$x \in S \cap H$, we have
$$\left\{ \begin{array}{l} x_1 = \psi(x_3, x_4)^2 +x_3^2 -x_4^2,\\
\mathbf{x} \in H. \\
\end{array} \right.$$
The classification of (real) quadratic surfaces implies that this
can be an elliptic paraboloid, a parabolic cylinder or a hyperbolic
paraboloid. An elliptic paraboloid contains no lines and the
corresponding case is trivial. If $S\cap H$ is a parabolic cylinder,
then all lines on it are parallel. It is straightforward that there
are $O(k^{2 + 2 \alpha})$ points in $P$ that lie on it (by counting
possible pairs $(x_3, x_4)$). Hence there are $O(k^{1 + 2\alpha})$
lines in $L$ that are fully contained in $S \bigcap H$. In the rest
of the discussion we assume $S \cap H$ is a hyperbolic
paraboloid.\ignore{ There are no lines on an elliptic paraboloid. If
$S \cap b$ is a parabolic cylinder then the points $x$ on $S \cap b$
satisfies that $x_2 = \psi(x_3, x_4)$ and $x_1 = \psi' (x_3, x_4)^2$
where $\psi'$ is some linear form. Note that in this case the
direction of lines on $S \cap H$ is fixed and each line in $L$
passes through $\Theta(k)$ points in $P$. So we only need to count
the number of points lies on $S \cap H$. This is obviously bounded
by the total number of possible $(x_3, x_4)$, which is $\Theta(k^{2
+ 2\alpha})$. Hence the total number of lines on $S \cap b$ is
$\lesssim k^{1 + 2 \alpha}$. Finally, if $S \cap b$ is a hyperbolic
paraboloid, then again} In this case, similarly to the case of the
one-sheeted hyperboloid, there are at most two lines with the same
direction. Moreover, the direction $(v_1, v_2, v_3, v_4)$ of any
line on $S \cap H$ satisfies $v_2 = \widetilde{\psi}(v_3, v_4)$ and
$v_4^2 = v_2^2 + v_3^2$ (where we let $\widetilde{\psi}$ denote the
linear homogeneous part of $\psi$). Thus once we fix $v_1$ and
``$v_3$ or $v_4$'' (depending on $\widetilde{\psi}$), we have
limited the possible direction $(v_1, v_2, v_3, v_4)$ in a set with
$\leq 2$ elements. Hence there are $O(k^{1 + 3 \alpha})$ lines that
are fully contained in $S \cap H$.

Finally, we show that for $\alpha<\frac 1 2$, the number of
incidences is (asymptotically) larger than
$\Theta\left(m^{2/5}n^{4/5}+m^{1/2}n^{1/2}q^{1/4}+m^{2/3}n^{1/3}+m+n\right)$,
which is the bound of Theorem~\ref{th:shso}, with $m=
\Theta(k^{3+3\alpha}), n=\Theta(k^{2+4\alpha}), q=O(k^{1+3\alpha}),$
and $s=O(1)$. We have
$$m^{2/5}n^{4/5}=O(k^{\frac {6+6\alpha+8+16\alpha} 5})=O(k^{\frac {14+22\alpha} 5}),$$
and the exponent is smaller than $3+4\alpha$, as $\alpha<\frac 1 2$.
Similarly,
$$m^{1/2}n^{1/2}q^{1/4}=O(k^{\frac {6+6\alpha+4+8\alpha+1+3\alpha} 4})=O(k^{\frac {11+17\alpha} 4}),$$
and the exponent is smaller than $3+4\alpha$, as $\alpha<\frac 1
2<1$. Similarly,
$$m^{2/3}n^{1/3}=O(k^{\frac {6+6\alpha+2+4\alpha} 3})=O(k^{\frac {8+10\alpha} 3}),$$
and the exponent is smaller than $3+4\alpha$ for every $\alpha$.
Since both $m$ and $n$ are $O(k^{3+4\alpha})$, the claim is proved.
\proofend

\section{Proof of Theorem~\ref{th:mainthree}}
The proof of Theorem~\ref{th:mainthree} follows easily by
Corollary~\ref{co:main}, together with the following lemma.

\begin {lemma}
\label{le:proj} Let $L$ be a set of $n$ lines in $\reals^4$ such
that at most $s$ lines lie on a common 2-flat. There exists a
projection from $\reals^4$ onto a hyperplane $H\subset \reals^4$,
such that at most $s$ lines lie on any common plane in $H$.
\end {lemma}

\noindent{\bf Proof of Lemma~\ref{le:proj}.} Let $\pi_1, \ldots,
\pi_k$ denote the set of 2-flats containing at least two lines in
$L$, then $k\le \binom n 2$. For a generic hyperplane $H \subset
\reals^4$, the projection $p : \reals^4 \to H$ maps $\pi_i$ onto a
plane $\pi'_i$ contained in $H$. We pick, as we may, a hyperplane
$H$, so that $p$ is bijective on $\pi_1,\ldots, \pi_k$. Denote by
$L'$ the set of projected lines in $\reals^3$. It is easy to verify
that the set of planes in $H$ containing at least two lines in $L'$
consists precisely of $\pi'_1,\ldots, \pi'_k$. Moreover, the number
of lines in $L'$ that are contained in $\pi'_i$ is equal to the
number of lines in $L$ that are contained in $\pi_i$, thus
completing the proof. \proofend

\section{Discussion and open questions}
In Corollary~\ref{co:main}, we show a concrete irreducible quadratic
hypersurface $S$ in $\reals^4$, together with a set of $m$ points
and $n$ lines that lie on $S$, for $n^{9/8}<m< n^{3/2}$, such that
(i) the number of lines in any common 2-flat is $O(1)$, (ii) the
number of lines in any common hyperplane is $O(n/m^{1/3})$, and
(iii) the number of incidences between the points and lines is
$\Omega(m^{2/3}n^{1/2})$, which is asymptotically larger than
$\Theta(m^{2/5}n^{4/5}+m^{1/2}n^{1/2}q^{1/4}+m^{2/3}n^{1/3}+m+n)$ in
this regime of $m$ and $n$. A natural question is to extend this
result to other regimes by a similar construction. The condition (i)
is natural and should not be hard to achieve, since if a plane is
not contained in a quadratic hypersurface, then by the generalized
version of B\'ezout's theorem~\cite{Fu84} it can contain at most two
lines. Here are a few natural questions that arise
\begin {enumerate}
\item
Can we generalize our construction, such that in (ii) we are allowed
to have a more general $q$, not necessarily $\sim n/m^{1/3}$, s.t.
the number of lines in any common hyperplane is $O(q)$, and we still
get a lower bound of incidences asymptotically larger than
$\Theta(m^{2/5}n^{4/5}+m^{1/2}n^{1/2}q^{1/4}+m^{2/3}n^{1/3}+m+n)$?
\item
Can we find a similar construction when $m<n^{9/8}$?
\item
How powerful is the natural generalization of this construction for $\reals^d$, when $d>4$?
Notice that for $d>4$, finding the precise bound for the number of
incidences between a set $P$ of $m$ points and a set $L$ of $n$
lines in $\reals^d$ is already an interesting open question.
It is probably too early for us to answer this question before we find the correct bound.
\item
In three dimensions, it remains a question to determine if Theorem~\ref{th:mainthree} is
tight.
\begin {question}
\label{conj} Let $P$ be a set of $m$ distinct points and $L$ a set
of $n$ distinct lines in $\reals^3$, and assume that no plane
contains more than $s=O(1)$ lines of $L$. Then what is a good or tight upper bound of $I(P, L)$? Would $O(m^{\frac{2}{3}}n^{\frac{1}{2}} + m + n)$ suffice?
\end {question}

We do not know the answer to this question yet. It seems to require new techniques.
\end {enumerate}


\end{document}